# Comparing different approaches to the concept of determinant: the multiplicative approach is the best!

**Peteris Daugulis**
*Daugavpils University, Latvia*

**Mikhail H Klin**
*Ben-Gurion University of the Negev, Israel*

**Anita Sondore**
*Daugavpils University, Latvia*

**Abstract.** *The determinant is traditionally introduced through permutation expansions, multilinearity, or recursive expansion formulas—approaches that, while rigorous, often obscure its conceptual significance. Instead, an approach is advocated wherein the determinant is defined via its natural characterization as a multiplicative map on matrices, reflecting its role as a structural invariant of linear transformations under composition. This perspective defers combinatorial machinery until after the core properties are established, thereby lowering the cognitive barrier for students encountering the concept for the first time. The equivalence of this characterization with classical formulations is shown, and it is argued that this framework better aligns with the structural emphasis of modern mathematical education.*

## Introduction

The matrix determinant is a central concept in linear algebra with far- reaching applications across mathematics and applications in physics, engineering, and data science. For any square matrix **A** of size n×n, the determinant is a scalar value, commonly denoted det(**A**), that encapsulates important algebraic and geometric properties of **A**. One of its fundamental properties is the characterization of invertibility: a matrix is invertible if and only if its determinant is nonzero; a vanishing determinant indicates singularity and the absence of an inverse matrix (Lax, 2007).

Geometrically, the determinant measures the scaling factor of n-dimensional volume induced by the linear transformation represented by **A**. In two dimensions, this corresponds to area scaling; in three dimensions, to volume scaling (Hannah, 1996).

Beyond invertibility and geometry, the determinant plays an important role in eigenvalue theory, where eigenvalues can be computed as the roots of a determinant- characteristic polynomial of the matrix (Axler, 2023; Hoffman & Kunze, 1971; Meyer, 2000; Zavalo et al., 1974).

Historically, the determinant emerged in the classical solution formulas for systems of linear equations, where it appears in the denominator and vanishes precisely when the system fails to have a unique solution. It is called the traditional combinatorial definition of the determinant (Freud, 2024). It is often the first tool students encounter for computing low-order determinants and is commonly presented as an introductory definition.

Subsequently the determinant was recognized to admit an axiomatic characterization aligned with its geometric meaning: it is the unique function on square matrices that is multilinear and alternating in the rows (or columns) and normalized so that the determinant of the identity matrix equals one. These properties determine the determinant uniquely (Artin, 1944; Artin, 1957). It is called the traditional multilinear definition of the determinant.

In addition to these structural characterizations, the determinant admits a recursive formulation via cofactor (Laplace) expansion along any row or column, providing an explicit computational definition. This definition can be justified by inductively solving systems of linear equations. It is called the traditional recursive definition of the determinant.

In mathematics, one may distinguish between computational definitions, which specify explicit methods of evaluation, and axiomatic definitions, which characterize an object through a system of defining properties. Note that the combinatorial and recursive definitions are computational definitions, while the multilinear definition is an axiomatic definition.

A central, well-known property of the determinant - and the primary focus of this article, is its multiplicativity:

$$\det(\mathbf{AB}) = \det(\mathbf{A}) \det(\mathbf{B}), \text{ for any } \mathbf{A}, \mathbf{B} \in Mat(n, k),$$

where *Mat(n, k)* is the set of n×n matrices over the ring *k*. This article develops a systematic exposition of a definition of the matrix determinant that takes this multiplicative property as a fundamental axiom. The underlying motivation for this approach is proposed, the classical properties of the determinant within this framework are established through rigorous proofs, and the conceptual and pedagogical advantages of this perspective are analyzed. Main results related to the message of this article can be found - both implicitly and explicitly - in the works cited below, as well as in numerous other sources. However, those presentations are often excessively sophisticated, and thus not well suited for university educators. The contribution of this article lies in providing a clear, self-contained exposition designed to be accessible to a broad audience of instructors teaching linear algebra at universities and colleges.

## 1. Traditional definitions and their criticism

### 1.1. The combinatorial definition

In mathematical modelling, observed data and unknown quantities are often related through first-degree polynomial relationships, giving rise to linear equations and systems of linear equations. By organizing data and unknowns into the natural data structures of linear algebra—vectors and matrices—these relationships can be represented compactly as systems of linear equations (SLEs). Linear functions constitute the simplest and most fundamental class of approximating functions, which largely explains the practical importance of SLEs across scientific and engineering applications. Indeed, systems of linear equations arise in EVERY branch of science, engineering, economics, and data analysis (Sahu, 2026). Consequently, the development and study of methods for solving such systems have been among the principal driving forces behind the emergence and evolution of linear algebra.

The general form of an SLE is defined as the system

$$\begin{cases} a_{11}x_1 + a_{12}x_2 + \cdots + a_{1n}x_n = b_1 \\ \cdots \\ a_{m1}x_1 + a_{m2}x_1 + \cdots + a_{mn}x_1 = b_m \end{cases},$$

where $x_j$– unknowns; $a_{ij}$ coefficients; $b_i$- constant terms. When solving square SLEs (m= n), it can be observed that the unknowns are rational functions of the system coefficients and constant terms.

Let us consider small- sized SLEs.

n=1. Let us consider the SLE

$$\{a_{11}x_1 = b_1.$$

It has a solution if $a_{11} \neq 0$; $x_1 = \frac{b_1}{a_{11}}$. One can define $\det[a_{11}] = a_{11}$. In these notations,

$$x_1 = \frac{\det[b_1]}{\det[a_{11}]}.$$

n=2. Let us solve the SLE **Ax=b**,

$$\begin{cases} a_{11}x_1 + a_{12}x_2 = b_1 \\ a_{21}x_1 + a_{22}x_2 = b_2 \end{cases}.$$

For simplicity assume the general position - divisions by all polynomials are possible. After some manipulations the solution is obtained:

$$x_1 = \frac{b_1a_{22} - b_2a_{12}}{a_{11}a_{22} - a_{12}a_{21}};\ x_2 = \frac{b_2a_{11} - b_1a_{21}}{a_{11}a_{22} - a_{12}a_{21}}.$$

Exactly one solution exists if $a_{11}a_{22} - a_{12}a_{21} \neq 0$. One can define

$$\det\begin{bmatrix} a_{11} & a_{12} \\ a_{21} & a_{22} \end{bmatrix} = a_{11}a_{22} - a_{12}a_{21}.$$

For example, $\det\begin{bmatrix} 1 & 2 \\ 3 & 4 \end{bmatrix} = 1\cdot4-2\cdot3=-2$.

In these notations one gets formulas (1),

$$x_1 = \frac{\det\begin{bmatrix} b_1 & a_{12} \\ b_2 & a_{22} \end{bmatrix}}{\det(\mathbf{A})};\ x_2 = \frac{\det\begin{bmatrix} a_{11} & b_1 \\ a_{21} & b_2 \end{bmatrix}}{\det(\mathbf{A})}. \quad (1)$$

Example. The solution of the system

$$\begin{cases} x_1 + 2x_2 = 5 \\ 3x_1 + 4x_2 = 6 \end{cases}$$

is

$$x_1 = \frac{\det\begin{bmatrix} 5 & 2 \\ 6 & 4 \end{bmatrix}}{\det\begin{bmatrix} 1 & 2 \\ 3 & 4 \end{bmatrix}} = -4,\ x_2 = \frac{\det\begin{bmatrix} 1 & 5 \\ 3 & 6 \end{bmatrix}}{\det\begin{bmatrix} 1 & 2 \\ 3 & 4 \end{bmatrix}} = \frac{9}{2} = 4.5.$$

n=3. To solve a 3×3 system of linear equations

**Ax=b**,

that is,

$$\begin{cases} a_{11}x_1 + a_{12}x_2 + a_{13}x_3 = b_1 \\ a_{21}x_1 + a_{22}x_2 + a_{23}x_3 = b_2 \\ a_{31}x_1 + a_{32}x_1 + a_{33}x_1 = b_3 \end{cases},$$

additional work is required. A natural inductive strategy is to use one of the equations to express one variable in terms of the others, and then substitute this expression into the remaining equations. This reduces the problem to solving a resulting 2×2 system of linear equations. Another classical approach is to use elementary row operations to eliminate variables and thereby reduce the system to a simpler form. The solution of such SLE can be expressed in the form

$$x_i = \frac{\det(\mathbf{A}_\mathrm{i})}{\det(\mathbf{A})},\ i \in \{1; 2; 3\},$$

where $\mathbf{A}_\text{i}$ is obtained from $\mathbf{A}$ by replacing the i-th column of $\mathbf{A}$ with $\mathbf{b}$, and for any matrix $\mathbf{A}$:

$$\det(\mathbf{A}) = a_{11}a_{22}a_{33} + a_{12}a_{23}a_{31} + a_{13}a_{21}a_{32} - a_{13}a_{22}a_{31} - a_{11}a_{23}a_{32} - a_{12}a_{21}a_{33}.$$

The following mnemonic rule – the Sarrus's rule, is taught to memorize the computation of determinants of 3×3 matrices, see Figure 1 (Lay et al., 2016). The first two columns of matrix $\mathbf{A}$ are duplicated and placed immediately to the right of the matrix. Products of triplets along the solid lines have a positive sign, while products of triplets along the dashed lines have a negative sign, see Figure 1. The determinant is equal to the sum of the six terms corresponding to these six triplets. No explanation of naturality or motivation for this rule is provided.

**Figure 1**

*A mnemonic picture to memorize the triple products*

$$\begin{bmatrix} a_{11} & a_{12} & a_{13} \\ a_{21} & a_{22} & a_{23} \\ a_{31} & a_{32} & a_{33} \end{bmatrix} \begin{matrix} a_{11} & a_{12} \\ a_{21} & a_{22} \\ a_{31} & a_{32} \end{matrix}$$

*Example.* Compute $\det \begin{bmatrix} 2 & -1 & 0 \\ 3 & 1 & 2 \\ 4 & 2 & -1 \end{bmatrix}$.

To apply Sarrus' rule,, the first two columns of determinant are copied directly to the right of the original, see Figure 2:

**Figure 2**

*Sarrus' rule lines for given determinant*

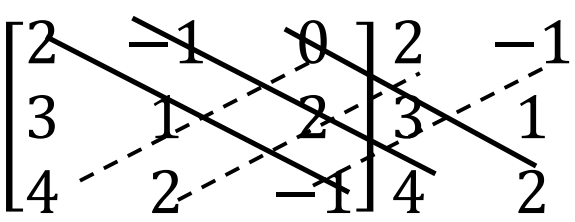

The determinant is then calculated by adding the products of the lines parallel to the main diagonal (solid lines) and subtracting the products of the lines paralles to the anti-diagonal (dashed lines), see Figure 2. One gets the answer

$$2·1·(-1)+(-1)\ 2·4+0·3·2-0·1·4-(-1)\ 3·(-1)-2·2·2=-21.$$

As the size of the system of linear equations increases, one can observe recurring patterns and use them to generalize its structure. This leads to the definition of a function

$$\det: Mat(n, k) \rightarrow k,$$

which is defined in a uniform way for all $n$. The classic combinatorial definition of determinant of a $n \times n$ matrix is given by the formula

$$\det(\mathbf{A}) = \sum_{\sigma\in\Sigma_n} \epsilon(\sigma) a_{1\sigma(1)}a_{2\sigma(2)}\cdots a_{n\sigma(n)} = \sum_{\sigma\in\Sigma_n} \epsilon(\sigma) a_{\sigma(1)1}a_{\sigma(2)2}\cdots a_{\sigma(n)n},$$

where the summation is over all permutations of the $n$-set, elements of $\Sigma_n$, and $\epsilon(\sigma)$ is the sign of a permutation $\sigma$ (Freud, 2024). Formulas of type (1) are used since the 18th century under the name of Cramer's rules (Lay et al., 2016).

The traditional combinatorial definition of the determinant provides an explicit closed-form formula and reveals important connections with discrete mathematics. From the theoretical perspective, it is elegant and complete, defining the determinant through a single expression. Nevertheless, the combinatorial definition of the determinant often strikes first-year undergraduates as unintuitive and opaque (Axler, 1995). By requiring a prerequisite mastery of abstract algebra- specifically the mechanics of permutations and their parity- this approach demands a background that is frequently absent from the early linear algebra curriculum. For students in non-mathematics majors, the combinatorial formula can be difficult to follow and understand. Its status as a standard introductory definition can be questioned because of its computational impracticality for large matrices and its weak algebraic motivation. This definition is practical only for 2×2 matrices. Furthermore, there exists a pedagogical dissonance between the relatively simple algebraic properties of the determinant and the complexity of its combinatorial origins.

### 1.2.The recursive definition

Solving SLEs inductively by using solution formulas for SLE of smaller sizes, it was observed that the determinants — polynomial expressions appearing in the numerators of the formulas for the unknowns — can themselves be expressed via determinants of matrices of smaller sizes. More precisely, the determinant of an n×n matrix can be written as a linear combination of determinants of (n−1)×(n−1) matrices obtained from the original matrix by deleting one row and one column (Chasnov, 2024; Lay et al., 2016).

Use the notation $\mathbf{A}_{ij}$ - the matrix obtained from **A** by deleting the i-th row and j-th column. Given $\mathbf{A}= [a_{ij}]_{n,n} \in Mat(n,k)$, a special case of such a linear combination is the expansion *along the first column:*

$$\det(\mathbf{A}) = \sum_{i=1}^{n} (-1)^{i+1} a_{i1} \det(\mathbf{A}_{i1}) .$$

*Example*. $\det\begin{bmatrix} a_{11} & a_{12} \\ a_{21} & a_{22} \end{bmatrix} = a_{11}\det[a_{22}] - a_{21}\det[a_{12}] = a_{11}a_{22} - a_{21}a_{12}$.

These expansion formulas allow one to define the determinant recursively- beginning with 1×1 matrices and proceeding inductively to matrices of arbitrary size. It is easily checked that in the cases of 2×2 and 3×3 matrices, this recursive construction reproduces the familiar combinatorial formulas for the determinant.

However, despite its computational usefulness and its historical origin in the theory of SLEs, the recursive definition largely obscures the conceptual meaning of the determinant. Within this approach, the determinant appears merely as a polynomial generated through repeated expansion formulas, with little indication of why such an expression should be regarded as natural or fundamental. In particular, the recursive construction conceals the determinant's essential algebraic role as the unique multiplicative scalar canonically associated with a linear transformation, as well as its geometric interpretation as the scaling factor of oriented volumes. As a result, students often perceive recursive determinant formulas as complicated combinatorial identities and mathematical tricks rather than as manifestations of a concept with clear algebraic and geometric significance.

### 1.3.The multilinear definition

Another definition of the determinant, now called the multilinear definition, emerged in the nineteenth century together with the development of linear algebra, multilinear algebra and geometry. It was observed in the early 19th century by Augustin-Louis Cauchy that the geometric volume of a parallelepiped is equal, up to a sign, to the determinant of a matrix whose columns are the vectors defining the solid's edges— a realization that provides a natural geometric definition of the determinant through its alternating multilinear properties. Rather than defining the determinant recursively through expansion formulas, this approach characterizes it through its geometric interpretation (Lax, 2007; Strang, 2016; Zavalo et al., 1974).

In the multilinear approach, the determinant is regarded as a function of the columns (or rows) of a matrix which is linear in each argument separately and changes sign whenever two rows are interchanged. These properties reflect the geometric intuition that the determinant measures oriented volume. Multilinearity describes how volumes behave under scaling and decomposition into parallelepipeds, while antisymmetry expresses the fact that exchanging two basis directions reverses orientation and therefore changes the sign of the volume. The vanishing of the determinant corresponds geometrically to the collapse of an n-dimensional parallelepiped into a lower-dimensional figure.

Let $\mathbf{A}\in$ *Mat(n, k)* be an n×n matrix and let $k^n$ be the set of columns of length $n$ having elements in the field $k$. Denote its columns by vectors $\mathbf{A}= [\mathbf{A}_1; \mathbf{A}_2;\ldots; \mathbf{A}_n]$, where each matrix $\mathbf{A}_j\in k^n$. According to the multilinear definition, the determinant is the function

$$\det: Mat(n, k) \rightarrow k,$$

that satisfies the following three properties:

1. **Multilinearity:** The function is linear with respect to each individual column vector. For any scalars $c_1, c_2\in k$ and vectors $\mathbf{v}_1, \mathbf{v}_2\in k^n$:

$$\det[\ldots,c_1 \mathbf{v}_1+ c_2 \mathbf{v}_2 ,\ldots]= c_1 \det[\ldots,\mathbf{v}_1 ,\ldots]+ c_2 \det[\ldots, \mathbf{v}_2 ,\ldots].$$

2. **Swapping Property:** If any two columns $\mathbf{v}$ of the matrix are identical, the determinant is zero:

$$\det[\ldots,\mathbf{v},\ldots, \mathbf{v} ,\ldots]=0.$$

This implies that swapping any two distinct columns changes the sign of the determinant:

$$\det[\ldots, \mathbf{A}_i; \ldots, \mathbf{A}_j;\ldots]=- \det[\ldots, \mathbf{A}_j; \ldots, \mathbf{A}_i;\ldots].$$

3. **Normalization:** The determinant of the identity matrix $\mathbf{E}_n$ is equal to 1:

$$\det(\mathbf{E}_n) =[ \mathbf{e}_1, \mathbf{e}_2 ,\ldots, \mathbf{e}_n]=1,$$

where $\mathbf{e}_j$ represents the j-th standard basis vector of $k^n$.

A fundamental theorem states that there exists a unique multilinear alternating function of the rows (equivalently, columns) of an $n\times n$ matrix that takes the value 1 on the identity matrix (Hoffman & Kunze, 1971). Thus, in this approach, the determinant is introduced not as a complicated polynomial formula, but as a function associated with the volume of the paralelepided determined by a sequence of vectors.

Nevertheless, despite its geometric elegance, the traditional multilinear definition has important pedagogical and structural shortcomings. It largely abandons the original motivation coming from systems of linear equations and replaces it with an axiomatic description whose naturality is not immediately evident to beginners. Important linear algebra techniques such as Cramer's rule are hidden. The determinant is introduced as a special function on ordered tuples of vectors, rather than as an intrinsic invariant of a linear operator itself. Consequently, the multilinear definition becomes closely tied to manipulations of rows, columns, and coordinates, placing disproportionate emphasis on geometric aspects, while the determinant's underlying linear-algebraic origin and conceptual significance remain largely concealed. This may be a learning obstruction for students in engineering or computer science, for whom the primary interaction with matrices involves arrays of numbers modelling discrete, interconnected physical variables (e.g., node voltages or structural load balances). Educational researchers note that this forces a "cognitive rupture," creating abstract barriers that hide the underlying linear equations behind unnecessary geometric intuition and machinery (Lin et al., 2024). In particular, the multiplicative nature of the determinant- arguably its most fundamental algebraic property- appears only later as a nontrivial theorem rather than as a guiding principle of the construction. As a result, students often gain geometric intuition about signed volumes without fully understanding why the determinant occupies such a central position in linear algebra.

## 2. The multiplicativity-based definition of the determinant

### 2.1. The multiplicative property

These pedagogical challenges underscore the value of an alternative approach: a definition that begins with a concise, technically modest algebraic description rooted in the basic linear algebra framework. Such definitions typically arise by characterizing a function as a morphism between appropriate algebraic structures.

In both purely mathematical and practical applications, as well as in mathematical modeling, it is useful to study functions from the considered complicated structures to sets of numbers. When studying algebraic structures such as matrices, functions that are morphisms from these structures to their underlying rings or fields are important. These morphisms transform algebraic structures into numbers and relate the operations of these structures to number operations. Such functions assign to a complex multidimensional structure a single number that condenses its properties.

A central property of the determinant is the multiplicative property

$$\det(\mathbf{AB}) = \det(\mathbf{A}) \det(\mathbf{B}), \text{ for any } \mathbf{A}, \mathbf{B} \in Mat(n, k),$$

where $Mat(n, k)$ is the set of $n \times n$ matrices over the ring $k$.

It is known (Brenner, 1968; Dieudonne, 1943; Vaccarino, 2009) that the classical determinant can be characterized as the unique function on square matrices satisfying multiplicativity together with suitable normalization conditions. In this formulation, multiplicativity is taken as the defining structural property, from which all the standard determinant definitions follow. It lifts the multiplicative property from an accidental numerical coincidence into a mandatory law. Despite the conceptual clarity of this approach, it is seldom developed in detail in standard textbooks and is largely absent from undergraduate linear algebra curricula. For a textbook discussion of a closely related definition, see (Kostrikin, 2000). It must be noted that the multiplicative definition of the determinant can be generalized over infinite-dimensional spaces (Gohberg et al., 2000).This feature is useful in studying areas such as mathematical physics and machine learning (Elizalde et al., 1994; Kulesza & Taskar, 2012).

### 2.2. Motivations from the process of solving systems of linear equations

Focusing on solving systems of linear equations and their interpretation in terms of matrices, the authors review considerations that can be utilized for defining a useful function

$$d: Mat(n, k) \to k.$$

#### 2.2.1 Denominators

By analogy with the equation

$$ax=b, \text{ where } a \neq 0,$$

with the solution $x = \frac{b}{a}$, the values of the function *d* must coincide with the denominators of the unknowns in the general case of the formulas, this also has nice traces from the historical approach:

$$\mathbf{Ax=b} \quad \Rightarrow \qquad \mathbf{x=A^{-1}\,b} = \frac{1}{d(\mathbf{A})}\begin{bmatrix} z_1 \\ z_2 \\ \cdots \\ z_n \end{bmatrix}.$$

One can think that $\frac{1}{d(\mathbf{A})}$ is a number that characterizes $\mathbf{A^{-1}}$.

### 2.2.2 Multiplicativity

Another motivation for the multiplicative property of the *d* function can be based on the SLE interpretation of square matrices (Daugulis, 2026).

A useful technique for solving systems of linear equations is a linear change of variables. To solve a system

$$\mathbf{Ax=b},$$

it may be advantageous to introduce a substitution of the form $\mathbf{x} = \mathbf{By}$, where $\mathbf{B}$ is an invertible matrix. This transforms the original system into

$$(\mathbf{AB})\mathbf{y} = \mathbf{b},$$

which is then solved with respect to $\mathbf{y}$, after which x is recovered via $\mathbf{x} = \mathbf{By}$. Viewing systems through such matrix compositions provides a natural motivation for the multiplicative property, as successive linear transformations correspond directly to products of matrices.

Consider the SLE

$$(\mathbf{AB})\mathbf{x} = \mathbf{b},$$

where $\mathbf{AB}$ is invertible. It can be solved in at least two different ways.

•In one step, then the denominator in the unknown formulas will be $d(\mathbf{AB})$:

$$\mathbf{x} = \frac{1}{d(\mathbf{AB})}\begin{bmatrix} z_1 \\ z_2 \\ \cdots \\ z_n \end{bmatrix}.$$

•In two steps:
the first step – one solves the equation

$$\mathbf{ABx=A(Bx)},$$

by substitution $\mathbf{x} \to \mathbf{y} = \mathbf{Bx}$. One gets the system

$$\begin{cases} \mathbf{Ay} = \mathbf{b} \\ \mathbf{y} = \mathbf{Bx} \end{cases},$$

solving the SLE with respect to $\mathbf{y}$, the denominator will be $d(\mathbf{A})$:

$$\mathbf{y} = \frac{1}{d(\mathbf{A})}\begin{bmatrix} y_1 \\ y_2 \\ \cdots \\ y_n \end{bmatrix},$$

the second step - solving the SLE with respect to **x**

$$\mathbf{Bx} = \mathbf{y} = \frac{1}{d(\mathbf{A})}\begin{bmatrix} y_1 \\ y_2 \\ \cdots \\ y_n \end{bmatrix},$$

one gets

$$\mathbf{x} = \frac{1}{d(\mathbf{B})} \cdot \frac{1}{d(\mathbf{A})}\begin{bmatrix} z_1 \\ z_2 \\ \cdots \\ z_n \end{bmatrix}.$$

Comparing denominators, the following relation is obtained

$$\frac{1}{d(\mathbf{AB})} = \frac{1}{d(\mathbf{A})} \cdot \frac{1}{d(\mathbf{B})}$$

from which follows the multiplicativity relation

$$d(\mathbf{AB}) = d(\mathbf{A}) \cdot d(\mathbf{B}).$$

**2.3. Consequences of multiplicativity for elementary matrices**

One can show that the multiplicative property implies certain values of the $d$ function on elementary matrices. Below we refer to the standard notation, related to the elementary row operations, see (Strang, 2016).

Properties of rings and solutions of functional equations are used.

Let $d$: $Mat(n, k) \to k$ be a morphism of multiplicative monoids, a nonconstant polynomial function of matrix elements. Then

1. $d(\mathbf{E}_n) = 1$, $d(\mathbf{0}_n) = 0$, where $\mathbf{0}_n$ is the zero matrix;
2. $\mathbf{M} \in Mat(n, k)$ is noninvertible iff $d(\mathbf{M}) = 0$;
3. for all $p \in \{1, ..., n\}$, $\lambda \in k$, there exists $m \in$ N such that

$$d(\mathbf{R}_p(\lambda)) = \lambda^m,$$

where $\mathbf{R}_p(\lambda)$ is the matrix of the elementary row operation of the second kind (ERO2);

4. for all $p, q \in \{1, ..., n\}$, $p \neq q$, for all $\lambda \in k$,

$$d(\mathbf{R}_{pq}(\lambda)) = 1,$$

where $\mathbf{R}_{pq}(\lambda)$ is the matrix of the elementary row operation of the third kind (ERO3);

5. if for all $p \in \{1, ..., n\}$, $\lambda \in k$, $d(\mathbf{R}_p(\lambda)) = \lambda^m$ with $m$ odd then for all $p, q$,

$$d(\mathbf{R}_{pq}) = -1,$$

where $\mathbf{R}_{pq}$ is the matrix of the elementary row operation of the first kind (ERO1).

**3. The definition**

One defines a determinant function using the multiplicativity property and its values on matrices of elementary row or column operations described in section 2.3. For all $n \in$ N the determinant function

$$\det : Mat(n, k) \to k$$

is defined according to the following rules:

1. for all $\mathbf{A}$, $\mathbf{B} \in Mat(n, k)$ $\det(\mathbf{AB}) = \det(\mathbf{A}) \det(\mathbf{B})$;
2. for all $p \in \{1, ..., n\}$, $\lambda \in k$ $\det(\mathbf{R}_p(\lambda)) = \lambda$.

Assuming these rules 1. and 2. statements given in section 2.3. implies the following properties 3.,4.,5.,6. which can be used in computations:

3. if $\mathbf{A}$ is noninvertible then $\det(\mathbf{A}) = 0$;
4. $\det(\mathbf{E}_n) = 1$;
5. for all $p, q \in \{1, ..., n\}$, $p \neq q$, $\det(\mathbf{R}_{pq}) = -1$;
6. for all $p, q$, $\forall \lambda \in k$, $\det(\mathbf{R}_{pq}(\lambda)) = 1$.

*Remark.* The rule 2. can be justified by the assumption that for a 1×1 matrix its determinant should be equal to the value of its sole element.

The value of $\det(\mathbf{A})$ can be calculated according to the following algorithm:

• if $\mathbf{A}$ is noninvertible then $\det(\mathbf{A}) = 0$, noninvertibility of $\mathbf{A}$ can be detected by transforming it to its row or column echelon form using elementary row or column operations;

• if $\mathbf{A}$ is invertible and $\mathbf{A} = \mathbf{P}_1 ... \mathbf{P}_l$, where each $\mathbf{P}_i$ is the matrix of an elementary row or column operation then

$$\det(\mathbf{A}) = \det(\mathbf{P}_1) ... \det(\mathbf{P}_l).$$

For such a definition to be justified and convenient to use, it is necessary to prove that the function is well- defined: if a matrix can be expressed as a product of elementary matrices in two different ways, then the determinant does not depend on it. It is also desirable to find an algorithm for calculating the determinant without using knowledge of matrix invertibility and factorization as a product of elementary matrices.

## 4. Consequences of the definition

Immediate consequences of the multiplicative definition of the determinant- such as changes induced by elementary row and column operations, as well as the determinants of diagonal and triangular matrices- are immediately noted.

It can be proved that the multilinearity of the determinant with respect to rows and columns follows from its multiplicativity. Further, using elementary row and column operations, one can give a slightly modified standard proof that the Laplace expansion is a consequence of row and column multilinearity- a property which can be derived from the multiplicative definition of the determinant. This argument provides a recursive definition of determinant in the classical sense, starting from 1x1 matrices. It proves that the multiplicative definition does not depend on matrix factorizations.

The classic combinatorial formula

$$\det(A) = \sum_{\sigma \in \Sigma_n} \epsilon(\sigma) a_{1\sigma(1)} a_{2\sigma(2)} \dots a_{n\sigma(n)} = \sum_{\sigma \in \Sigma_n} \epsilon(\sigma) a_{\sigma(1)1} a_{\sigma(2)2} \dots a_{\sigma(n)n},$$

where $\epsilon(\sigma)$ is the sign of $\sigma$, can be proved using the multiplicative definition by proving that the combinatorial formula holds for elementary matrices and the multiplicative property holds for multiplication by elementary matrices.

## Personal teaching reminiscences

The teaching experience of Mikhail Klin started in the 1970s. During the period 1978 - 1983 Mikhail Klin was working at Kaluga State Pedagogical Institute, named after K.E.Tsiolkovsky. employed at the Department of Algebra and Geometry, moving from the position of Senior Lecturer to Associate Professor.

The standards of teaching were high, and colleagues were receptive to innovation. Algebra and Number Theory was a four-semester course. The main textbook was the Russian edition (Zavalo et al, 1974) of the original Ukrainian text. The group of students consisted of a few dozens of young persons, from diverse countries, among them: Cuba, Ethiopia, Guinea-Bissau. A handful of students were from the Chechen Republic. Another half was from diverse places of Kaluga district.

Practically all students succeeded in learning the concept of a determinant at an axiomatic level. For many of them, it served as a window into the world of mathematical abstraction.

During the 1990s and subsequent years, Mikhail Klin taught a number of mathematics courses at Ben-Gurion University of the Negev in Beer Sheva, Israel. Notably, the textbook (Hoffman & Kunze, 1971) remains a widely used and highly regarded text in linear algebra at the university.

Decades of experience in teaching linear algebra and its many applications have led to numerous pedagogical observations and suggestions. In particular, there is room for innovation in the definitions and sequencing of some fundamental concepts. One such concept is the determinant of a matrix.

During the preparation of a recent textbook (Daugulis, 2026), it became apparent that some of the innovative aspects of the multiplicative property of the determinant are still unfamiliar to a portion of the mathematics education community in Latvia. This observation motivated us to present the underlying ideas and their pedagogical implications in greater detail.

## Discussion and conclusions

Authors have presented a definition of the determinant of a matrix that takes multiplicativity as the fundamental axiom. It is noted that definitions of this type have been treated in the literature as optional material, typically presented at the end of the corresponding sections. Our approach differs slightly from the similar definition given in (Kostrikin, 2000) by not requiring antisymmetry with respect to rows and columns. In many standard textbook treatments it is not emphasized that the antisymmetry follows as a consequence of multiplicativity.

In contemporary mathematical practice, determinants are rarely computed using explicit permutation expansions or recursive cofactor formulas. Instead, practical computations rely almost exclusively on elimination procedures, matrix factorizations, and numerical algorithms implemented in computer algebra systems and scientific software. At the same time, recent trends in mathematical education increasingly emphasize structural understanding, conceptual reasoning, and the interpretation of mathematical objects through their invariant properties rather than through lengthy symbolic manipulations. In this context, the traditional combinatorial introduction of the determinant may appear pedagogically unsatisfactory, since it often presents the determinant first as a complicated formal expression whose conceptual meaning emerges only later. The multiplicative viewpoint offers an alternative organization of the subject.

A definition based on multiplicativity has clear linear–algebraic motivations. The basic properties of the determinant can be established within one or two lectures. It is suggested that the determinant section be positioned immediately after the matrix representation of systems of linear equations is introduced.

Such an approach can be useful in teaching linear algebra and its applications, since many important results involving determinants rely precisely on their multiplicative property. It is logically and pedagogically more justified to derive the complicated combinatorial formula for the determinant from its simple multiplicative property, rather than the other way around. The multiplicative definition of the determinant shifts the learning focus from cumbersome computations to understanding the composition of linear mappings. Rather than viewing the determinant as a static sum of products, students would perceive it as a dynamic scaling factor that characterizes a linear operator. Treating the matrix as a whole, and the determinant as a function of it, helps students understand key properties like invertibility without getting lost in row-and-column details. For 2×2 matrices the combinatorial and recursive formulas can be derived in an ad hoc style. Computing determinants using triangulation should be the main general method for matrices of starting from the size 3×3.

The multiplicative definition (unlike the combinatorial or recursive definitions) can play an important role in functional analysis and mathematical physics, where variants of determinants for the infinite-dimensional case.

A suggestion for future work is to give an undergraduate-accessible proof of the well-definedness of the determinant that does not rely on Laplace expansion or recursive definitions.

## Acknowledgements

Partially supported by The Center for Advanced Studies in Mathematics, Ben-Gurion University of the Negev, Beer Sheva, Israel.

## Declaration of Generative AI and AI-assisted Technologies

During the preparation of this work, the authors used ChatGPT to assist with translation and language improvement. After using this tool, the authors reviewed and edited the content as needed and take full responsibility for the content of the publication.